\documentclass[11pt,letterpaper]{amsart}

\usepackage{graphicx}
\usepackage{mathtools,amsmath,amsthm,amsfonts,amssymb}
\usepackage{enumitem}
\usepackage{tikz}
\usetikzlibrary{patterns}
\usepackage{comment}

\usepackage[colorlinks,allcolors=blue,backref=page]{hyperref}
\usepackage{aliascnt}
\usepackage[capitalize]{cleveref}
\newcommand{\sbar}[1]{\overline{\vphantom{1}#1}}

\newtheorem{theorem}{Theorem}[section]

\newaliascnt{lemma}{theorem}
\newtheorem{lemma}[lemma]{Lemma}
\aliascntresetthe{lemma}

\newaliascnt{deftheorem}{theorem}

\aliascntresetthe{deftheorem}

\newaliascnt{proposition}{theorem}
\newtheorem{proposition}[proposition]{Proposition}
\aliascntresetthe{proposition}

\newaliascnt{corollary}{theorem}
\newtheorem{corollary}[corollary]{Corollary}
\aliascntresetthe{corollary}

\newaliascnt{claim}{theorem}

\aliascntresetthe{claim}

\newaliascnt{conjecture}{theorem}

\aliascntresetthe{conjecture}

\theoremstyle{definition}

\newaliascnt{definition}{theorem}
\newtheorem{definition}[definition]{Definition}
\aliascntresetthe{definition}

\newaliascnt{notation}{theorem}

\aliascntresetthe{notation}

\newaliascnt{example}{theorem}
\newtheorem{example}[example]{Example}
\aliascntresetthe{example}

\theoremstyle{remark}

\newaliascnt{remark}{theorem}
\newtheorem{remark}[remark]{Remark}
\aliascntresetthe{remark}

\newaliascnt{question}{theorem}

\aliascntresetthe{question}

\newaliascnt{fact}{theorem}

\aliascntresetthe{fact}

\renewcommand*{\backrefalt}[4]{%
  \ifcase #1 %
  \or
    \(\uparrow\)~#2%
  \else
    \(\uparrow\)~#2%
  \fi}

\usepackage{xcolor}
\definecolor{airforceblue}{rgb}{0.36, 0.54, 0.66}

\definecolor{zsxpurple}{RGB}{128,0,128}

\definecolor{Lightred}{RGB}{239,83,80}

\definecolor{a}{rgb}{0.5,0.75,0.69}

\title[On the asymmetry polynomial of Schubert varieties]
{On the asymmetry polynomial of Schubert varieties}

\author[G. Burrull]{Gaston Burrull}
\address{International Institute for Mathematical Sciences, Jilin University, Qianjin Avenue 2699, High-Tech Development Zone, Changchun, Jilin 130012, China}
\email{grburrul@gmail.com}

\author[Z. A. Li]{Zhitong Alex Li}
\address{Department of Mathematics, Dornsife College of Letters, Arts and Sciences, University of Southern California, Los Angeles, CA 90089, USA}
\email{ali24951@usc.edu}

\author[S. Zeng]{Shixuan Zeng}
\address{School of the Gifted Young, University of Science and Technology of China, No. 96 Jinzhai Road, Baohe District, Hefei, Anhui 230026, China}
\email{zsx0515@mail.ustc.edu.cn}

\author[X. Zhang]{Xinyuan Zhang}
\address{Department of Mathematics, Southern University of Science and Technology, 1088 Xueyuan Avenue, Shenzhen, Guangdong 518055, China}
\email{12312009@mail.sustech.edu.cn}

\date{}



\hypersetup{
  pdftitle={On the asymmetry polynomial of Schubert varieties},
  pdfauthor={Gaston Burrull, Zhitong Alex Li, Shixuan Zeng, and Xinyuan Zhang}
}

\begin{document}

\begin{abstract}
    Let \(W\) be a finite Weyl group with length function \(\ell\) and longest element \(w_0\).
    For \(w\in W\), let \(P_w(q)\) be the Poincar\'e polynomial of the Schubert variety \(X_w\).
    We introduce its \emph{asymmetry polynomial}
    \[
        D_w(q):=q^{\ell(w)}P_w(q^{-1})-P_w(q).
    \]

    Let \(n\geq4\) and \(W=S_n\) or \(W(B_n)\).
    For \(w\in W\) whose Kazhdan--Lusztig polynomial \(P_{\mathrm{id},w}(q)\) satisfies \(P_{\mathrm{id},w}(1)=2\), we prove the sharp bound
    \[
        [q^k]D_w(q)\leq[q^{k-1}]P_{w_0}(q)
        \quad\text{for}\quad 1\leq k<\ell(w)/2.
    \]

    For every \(n\geq3\) and \(1\leq h\leq2n-4\), we construct an element \(b=b_h^n\in W(B_n)\) with \(P_{\mathrm{id},b}(q)=1+q^h\), and determine the asymmetry polynomial \(D_b(q)\) and all Kazhdan--Lusztig polynomials \(P_{x,b}(q)\) for \(x\leq b\).
\end{abstract}
\maketitle
\section{Introduction}\label{section:introduction}

Let \(W\) be a finite Weyl group with Bruhat order \(\leq\),
identity \(\mathrm{id}\), length function \(\ell\), and
longest element \(w_0\).
Write \(X_w\) for the corresponding Schubert variety and
\(P_{x,w}(q)\) for the Kazhdan--Lusztig polynomial, for \(x\leq w\).
Write
\[
    P_w(q)=\sum_{x\leq w}q^{\ell(x)}=\sum_i a_iq^i
\]
for the ordinary Poincar\'e polynomial of \(X_w\), where
\(a_i=\dim H^{2i}(X_w,\mathbb Q)\).
By the Carrell--Peterson criterion, \(P_w(q)\) is palindromic
exactly when \(X_w\) is rationally smooth
\cite[Theorem~D]{BE09}.
To measure its failure of palindromicity, we introduce the \emph{asymmetry polynomial}
\[
    D_w(q)=q^{\ell(w)}P_w(q^{-1})-P_w(q).
\]

Write \([q^i]f(q)\) for the coefficient of \(q^i\) in a polynomial \(f(q)\).
By Bj\"orner and Ekedahl~\cite[Theorem~A]{BE09},
\[
    [q^i]D_w(q)=a_{\ell(w)-i}-a_i\geq0
    \qquad(0\leq i<\ell(w)/2).
\]
Moreover, when \(P_{\mathrm{id},w}(q)\neq1\), the first positive coefficient of \(D_w(q)\) occurs in the smallest positive degree of \(P_{\mathrm{id},w}(q)\) \cite[Theorem~D]{BE09}.

We focus on \(P_{\mathrm{id},w}(1)=2\), the smallest possible value for a rationally singular Schubert variety.
By nonnegativity of Kazhdan--Lusztig coefficients, this condition implies \(P_{\mathrm{id},w}(q)=1+q^h\) for some \(h\geq1\).
The monotonicity theorem gives \(P_{x,w}(q)\in\{1,1+q^h\}\) for every \(x\leq w\).
In type \(A\), Woo characterized these permutations by the irreducibility of their singular locus together with avoidance of six specified patterns \cite{Woo09}.
We study how large the asymmetry coefficients \([q^k]D_w(q)\) can be for \(1\leq k<\ell(w)/2\).

Our first result establishes sharp bounds in the symmetric group \(S_n=W(A_{n-1})\) and the signed permutation group \(W(B_n)\).

\begin{theorem}\label{thm:sharp-asymmetry}
    Let \(n\geq4\) and \(W=S_n\) or \(W(B_n)\).
    Suppose \(w\in W\) satisfies \(P_{\mathrm{id},w}(1)=2\).
    For \(1\leq k<\ell(w)/2\), we have the sharp bound
    \[
        [q^k]D_w(q)\leq[q^{k-1}]P_{w_0}(q).
    \]
    For each \(1\leq k\leq\lfloor(n-2)/2\rfloor\) in \(S_n\), and \(1\leq k\leq n-2\) in \(W(B_n)\), there is such an element \(w\) with \(k<\ell(w)/2\) attaining equality at degree \(k\).
\end{theorem}

Under the same hypothesis \(P_{\mathrm{id},w}(1)=2\), the
inequality holds in every finite Weyl group
(\Cref{cor:universal-asymmetry-bound}).
In Section~\ref{section:asymmetry}, we establish sharpness in types \(A\) and \(B\) via explicit families that attain equality in several consecutive degrees, and we compute their asymmetry polynomials (\Cref{thm:type-A-uniform,thm:type-B-sharp-family}).

\begin{example}\label{exa:Dw-in-S4}
    In one-line notation, let \(w=4231\in S_4\), so that \(\ell(w)=5\) and \(w_0=4321\).
    We have \(P_{\mathrm{id},w}(q)=1+q\), and
    \[
        \begin{aligned}
            P_{w_0}(q)
            &=1+3q+5q^2+6q^3+5q^4+3q^5+q^6,\\
            D_w(q)
            &=q+q^2-q^3-q^4.
        \end{aligned}
    \]
    Thus the bound in \Cref{thm:sharp-asymmetry} is attained
    in degree~\(1\) and is strict in degree~\(2\):
    \[
        [q]D_w(q)=1=[q^0]P_{w_0}(q),
        \qquad
        [q^2]D_w(q)=1<3=[q]P_{w_0}(q).
    \]
\end{example}
\begin{example}
    The hypothesis \(P_{\mathrm{id},w}(1)=2\) matters already in degree
    one. For \(z=34512\in S_5\), we have \(P_{\mathrm{id},z}(q)=1+2q\),
    but
    \[
        [q]D_z(q)=2>1=[q^0]P_{54321}(q).
    \]
\end{example}

We also determine the maximum number of positive coefficients of \(D_w(q)\) below the middle degree in types \(A\) and \(B\).
To state this, define
\[
    N(w)=\#\{i\in\mathbb Z:1\leq i<\ell(w)/2,\ [q^i]D_w(q)>0\}.
\]
Let
\[
    u_n=
    \begin{cases}
        n\,(n-2)(n-1)(n-3)\cdots1,&W=S_n,\\
        \overline1\cdots\overline{n-3}\,
        \overline{n-1}\,\overline{n-2}\,\overline n,
        &W=W(B_n).
    \end{cases}
\]
Here \(\overline j=-j\) denotes a signed entry, and empty strings are omitted.

\begin{corollary}\label{thm:number-asymmetries}
    Let \(n\geq4\), and let \(W=S_n\) or \(W=W(B_n)\).
    We have
    \(P_{\mathrm{id},u_n}(q)=1+q\) and
    \[
        N(u_n)=\max_{w\in W}N(w)
        =\left\lfloor\frac{\ell(w_0)}2\right\rfloor-1.
    \]
\end{corollary}


Finally, we study the first degree in which palindromicity fails.
Gaetz and Gao determined this degree for every singular Schubert variety in type \(A\) \cite[Theorem~1.6]{GG24}.
Our final result constructs elements in \(W(B_n)\) with \(P_{\mathrm{id},w}(q)=1+q^h\) for every \(1\leq h\leq2n-4\), and determines all their Kazhdan--Lusztig polynomials and asymmetry polynomials.
For \(m\geq1\), write \([m]_q=1+q+\cdots+q^{m-1}\).

\begin{theorem}\label{thm:type-B-height-spectrum}
    Let \(n\geq3\) and \(1\leq h\leq2n-4\).
    There are \(a_h^n\leq b_h^n\) in \(W(B_n)\) such that, for every \(x\leq b_h^n\),
    \[
         P_{x,b_h^n}(q)=
         \begin{cases}
             1+q^h,&x\leq a_h^n,\\
             1,&x\not\leq a_h^n.
         \end{cases}
    \]
    The rationally singular locus is the rationally smooth variety \(X_{a_h^n}\), and
    \[
        D_{b_h^n}(q)=(q^h-q^{2n-3})\prod_{j=1}^{n-2}[2j]_q.
    \]
\end{theorem}

We give the explicit construction and proof in Section~\ref{sec:type-B-families}, where we also identify the positive coefficients of \(D_{b_h^n}\).

For comparison, in type \(A_n\), whose Weyl group is \(S_{n+1}\),
the equality \(P_{\mathrm{id},w}(q)=1+q^h\) forces \(h\leq n-2\)
by \cite[Theorem~1.5]{GG24}. Thus for \(n-1\leq h\leq 2n-4\),
\Cref{thm:type-B-height-spectrum} realizes the Kazhdan--Lusztig polynomial \(1+q^h\) at the
identity in type \(B_n\), beyond the range possible at the identity
in type \(A_n\).
\subsection*{Acknowledgments}
We thank the organizers of the PACE 2026 program at BICMR, Peking University, where this project originated, for their hospitality and support.
\subsection*{AI disclosure}
The authors used ChatGPT (GPT-5.6, GPT-6, and GPT-6 Astra) and Yuanbao (Hy3 and Hy4) for literature searches, mathematical exploration, computation, writing, and review.
ChatGPT assisted with SageMath code and computational checks of Kazhdan--Lusztig polynomials, asymmetry formulas, and examples illustrating the main results.
For \Cref{thm:type-B-sharp-family}, ChatGPT helped work out the details of an author's incidence-variety construction; another author later identified a simplification and used GPT-6 Astra to help revise the proof.
ChatGPT was also used to review the geometric proof of \Cref{thm:type-B-height-spectrum}, including the flag projection, its fibers, and the transverse quadratic cone.
The authors checked and revised the mathematical arguments and take full responsibility for the final manuscript.

\section{Preliminaries}\label{section:preliminaries}

\subsection{Schubert varieties and Kazhdan--Lusztig polynomials}
\label{subsec:KL-background}

Let \(S\) denote the set of simple reflections of \(W\), and put
\([u,w]=\{x\in W:u\leq x\leq w\}\).
Choose a connected complex reductive group \(G\) with Weyl group
\(W\), a Borel subgroup \(B\), and a maximal torus
\(\mathcal T\subseteq B\). For a representative
\(\dot w\in N_G(\mathcal T)\), set
\[
     X_w^\circ=B\dot wB/B,
     \qquad X_w=\overline{X_w^\circ}.
\]
The cell \(X_w^\circ\) is an affine space of dimension \(\ell(w)\),
and
\[
    X_w=\bigsqcup_{x\leq w}X_x^\circ.
\]
Thus the cell decomposition gives
\(P_w(q)=\sum_{x\leq w}q^{\ell(x)}\).

The Kazhdan--Lusztig polynomial \(P_{x,w}(q)\) records the
local intersection cohomology of \(X_w\) along the Schubert cell
\(X_x^\circ\), with cohomological degrees divided by two
\cite{KL80}. We use the normalization in which a smooth point
contributes \(1\).
We have
\[
     P_{w,w}=1,\qquad P_{x,w}(0)=1,
     \qquad \deg P_{x,w}<\frac{\ell(w)-\ell(x)}2\quad(x<w),  
\]
by the defining properties of KL polynomials \cite[\S2]{KL79}.
The intersection-cohomology interpretation also gives nonnegative
integer coefficients.

We write \(\operatorname{Gr}(p,V)\) for the Grassmannian of
\(p\)-dimensional subspaces of \(V\),
\(\operatorname{Gr}(p,n)=\operatorname{Gr}(p,\mathbb C^n)\),
and \(\mathbb P(V)\) for the space of lines in \(V\).
A complete flag \(E_\bullet\) in an \(n\)-dimensional space has
\(\dim E_i=i\), including \(E_0=0\) and \(E_n=V\).

We write \(\operatorname{Sing}(X_w)\) for the singular locus of \(X_w\).

The variety \(X_w\) is rationally smooth along \(X_x^\circ\)
if and only if \(P_{x,w}(q)=1\). Thus \(X_w\) is rationally
smooth if this holds for every \(x\leq w\).

The combinatorial definitions of \(P_w\) and \(D_w\) also apply
to arbitrary Coxeter systems, since lower Bruhat intervals
are finite. We use this generality in
\Cref{prop:whole-locus-transfer}.

\begin{proposition}
\label{prop:monotonicity}
    Let \((W,S)\) be a Coxeter system and let \(u\leq v\leq w\). Then
    \[
        P_{u,w}(q)-P_{v,w}(q)\in\mathbb Z_{\geq0}[q].
    \]
    If \(P_{u,w}(1)=2\), then, for some integer \(r>0\),
    \[
        P_{u,w}(q)=1+q^r,
        \qquad P_{v,w}(q)\in\{1,1+q^r\}.
    \]
\end{proposition}

\begin{proof}
    The coefficientwise inequality holds for arbitrary Coxeter
    systems by \cite[\S5]{Pla17}; the finite Weyl group case
    is \cite[Corollary~4]{Ir88}.
    Taking \(v=w\) also gives positivity.
    If \(P_{u,w}(1)=2\), positivity and the constant term one give
    \(P_{u,w}=1+q^r\) for some \(r>0\). Monotonicity then leaves
    only the two stated possibilities for \(P_{v,w}\).
\end{proof}

\begin{definition}\label{def:ratsing}
    Let \((W,S)\) be a Coxeter system and \(w\in W\). Define
    \[
        \operatorname{RatSing}(w)
        :=
        \{x\leq w:P_{x,w}(q)\neq1\},
    \]
    and
    \[
        \operatorname{MaxRatSing}(w)
        :=
        \operatorname{Max}\bigl(\operatorname{RatSing}(w)\bigr).
    \]
    Here maximality is with respect to the Bruhat order.
\end{definition}

For a finite Weyl group, the rationally singular locus of \(X_w\) is
\[
    \operatorname{RatSing}(X_w)
    :=\bigcup_{x\in\operatorname{RatSing}(w)}X_x^\circ.
\]
Its irreducible components are the varieties \(X_v\) with \(v\in\operatorname{MaxRatSing}(w)\).

We compose permutations and signed permutations by
\((uv)(i)=u(v(i))\), and use one-line notation
\(y=y_1\cdots y_n\), where \(y_i=y(i)\).
We realize \(W(B_n)\) as signed permutations whose entries
have absolute values \(1,\ldots,n\), extended by \(y(-i)=-y(i)\),
and write \(\overline a=-a\).
In \(S_n\), the simple reflections are \(s_i=(i,i+1)\),
for \(1\leq i<n\). In \(W(B_n)\), right multiplication by
\(s_0\) changes the sign of the first entry, and right
multiplication by \(s_i\), for \(1\leq i<n\), exchanges
positions \(i,i+1\).

By the first part of \Cref{prop:monotonicity}, \(\operatorname{RatSing}(w)\) is a lower order ideal.

\begin{example}\label{ex:rational-vs-ordinary-singular}
    Let
    \[
        w=1\overline2\overline3 4,
        \qquad v=3214,
        \qquad u=21\overline3 4
        \qquad\text{in }W(B_4).
    \]
    For every \(x\leq w\),
    \[
        P_{x,w}(q)=
        \begin{cases}
            1+q^2,&x\leq v,\\
            1,&x\not\leq v.
        \end{cases}
    \]
    Consequently, \(\operatorname{RatSing}(w)=[\mathrm{id},v]\) and \(\operatorname{MaxRatSing}(w)=\{v\}\).

    Since all three signed permutations fix \(4\), they lie in the
    standard parabolic subgroup \(W(B_3)\). To compute the ordinary
    singular locus, we use the symplectic realization of this Coxeter
    system, of type \(C_3\), and project to the Lagrangian Grassmannian
    of a six-dimensional symplectic space.
    This Grassmannian parametrizes the three-dimensional
    subspaces on which the symplectic form vanishes.
    The image of \(X_w\) is the Schubert divisor
    \[
        Y=\{L:\dim(L\cap E_3)\geq1\},
    \]
    where \(E_3\) is the fixed Lagrangian subspace.
    In symmetric-matrix coordinates, this condition is \(\det A=0\), with \(A\) a symmetric \(3\times3\) matrix.
    Its singular locus is
    \[
        \operatorname{Sing}(Y)
        =\{L:\dim(L\cap E_3)\geq2\}.
    \]
    Since \(X_w\) is the inverse image of \(Y\) under the smooth flag projection, its singular locus is the inverse image of \(\operatorname{Sing}(Y)\), namely \(X_u\).
    Thus,
    \[
        \operatorname{RatSing}(X_w)=X_v
        \subsetneq X_u=\operatorname{Sing}(X_w).
    \]
    Indeed, \(P_{u,w}(q)=1\), so \(X_w\) is singular but rationally smooth along \(X_u^\circ\).
\end{example}

\begin{definition}
    For \(x\leq y\) with \(P_{x,y}(q)\neq 1\), we define the
    \emph{Kazhdan--Lusztig height} of the pair \((x,y)\) by
    \[
        h_{\mathrm{KL}}(x,y)
        :=
        \min\{j>0:[q^j]P_{x,y}(q)>0\}.
    \]
\end{definition}

\subsection{Upward degrees in the Bruhat graph}
\label{subsec:upward-degree-bruhat-graph}

For \(y\in W(B_n)\), the length is given by \cite[Proposition~8.1.1]{BB05}:
\[
    \ell(y)
    =
    \#\{i:y_i<0\}
    +\#\{i<j:y_i>y_j\}
    +\#\{i<j:-y_i>y_j\}.
\]
Recall that the set of reflections is
\[
    T=\{wsw^{-1}:w\in W,\ s\in S\}.
\]
In \(W(B_n)\), the \(n^2\) reflections act by right
multiplication as a sign change, an exchange of two entries,
or an exchange of two entries together with changes of both
signs \cite[Proposition~8.1.5]{BB05}.

The \emph{Bruhat graph} of \([\mathrm{id},y]\) has vertex set \([\mathrm{id},y]\), with an edge between \(v\) and \(tv\) whenever \(t\in T\) and \(v<tv\leq y\). These edges need not be covers in the Bruhat order.

\begin{definition}\label{def:upward-degree}
    For \(v\leq y\), set
    \[
        d_y(v):=\ell(y)-\ell(v)
        =\operatorname{codim}_{X_y}X_v
    \]
    and define the \emph{upward degree} of \(v\) in the Bruhat graph by
    \[
        \deg_y^+(v)
        :=\#\{t\in T:v<tv\leq y\}.
    \]
\end{definition}

The same neighbors are obtained by right multiplication, since
conjugation by \(v\) permutes the set of reflections.
By Deodhar's inequality (see \cite[Section~6.2]{BL00}),
\[
    \deg_y^+(v)\geq d_y(v).
\]
Moreover, the Carrell--Peterson criterion implies that \(X_y\) is rationally smooth if and only if equality holds for every \(v\leq y\).

The following identity relates Kazhdan--Lusztig polynomials to a weighted sum over upward edges in the Bruhat graph.
In particular, it gives a formula for the upward degree at maximal rationally singular elements, which we use to study equality in the coefficient bound in Section~\ref{subsec:sharp-bounds}.
\begin{lemma}[Kobayashi {\protect\cite[Lemma~8.1]{Kob13}}]
\label{lem:KL-derivative-identity}
    For \(x\leq y\),
    \[
        d_y(x)P_{x,y}(1)-2P'_{x,y}(1)
        =
        \sum_{\substack{t\in T\\x<tx\leq y}}P_{tx,y}(1).
    \]
\end{lemma}

If \(x\in\operatorname{MaxRatSing}(y)\), then \(P_{tx,y}(q)=1\) for every \(t\in T\) with \(x<tx\leq y\). Thus the right-hand side equals \(\deg_y^+(x)\).

\begin{example}
    For the pair \(v=3214\) and \(w=1\overline2\overline3 4\) from Example~\ref{ex:rational-vs-ordinary-singular}, we have \(\ell(w)=8\), \(\ell(v)=3\), and \(P_{v,w}(q)=1+q^2\).
    Since \(v\in\operatorname{MaxRatSing}(w)\), \Cref{lem:KL-derivative-identity} gives
    \[
        \deg_w^+(v)=5\cdot2-2\cdot2=6>d_w(v)=5.
    \]
\end{example}

\subsection{Polynomial notation and unimodality}
\label{subsec:polynomial-tools}

For the \(q\)-integers introduced above, define
\[
    [m]_q!=\prod_{j=1}^m[j]_q,
    \qquad [0]_q!=1.
\]
For \(0\leq b\leq a\), the
\emph{Gaussian binomial coefficient} is
\[
    \genfrac{[}{]}{0pt}{}{a}{b}_q
    =\frac{[a]_q!}{[b]_q![a-b]_q!}.
\]
The Schubert cell decompositions give the Poincar\'e polynomials
\([a]_q!\) for the complete flag variety of \(\mathbb C^a\)
and \(\genfrac{[}{]}{0pt}{}{a}{b}_q\) for \(\operatorname{Gr}(b,a)\).
For the longest elements in types \(A\) and \(B\),
\[
    \begin{aligned}
        W=S_n:\qquad
            &P_{w_0}(q)=[n]_q!,
            &\ell(w_0)&=\binom{n}{2},\\
        W=W(B_n):\qquad
            &P_{w_0}(q)=\prod_{j=1}^n[2j]_q,
            &\ell(w_0)&=n^2.
    \end{aligned}
\]

The Poincar\'e polynomial of a rationally smooth Schubert variety
is palindromic and unimodal \cite[Theorems~A and~D]{BE09}.
To prove positivity below the middle for our asymmetry
polynomials, we need the following stronger property.

\begin{definition}\label{def:strict-unimodality}
    A palindromic polynomial \(f(q)=\sum_{i=0}^d b_iq^i\) is
    \emph{strictly unimodal} if
    \[
         0<b_0<b_1<\cdots<b_{\lfloor d/2\rfloor}.
    \]
\end{definition}

\begin{lemma}\label{lem:multiply}
    Let \(f(q)\) be strictly unimodal and palindromic of degree \(d\).
    If \(1\leq m<d\), then \(f(q)[m+1]_q\) has the same properties.
    Moreover, if \(1\leq t\leq d\), then
    \[
         [q^j](1-q^t)f(q)>0
         \qquad\left(0\leq j<\frac{d+t}{2}\right).
    \]
\end{lemma}

\begin{proof}
    Extend \(b_i\) by zero outside \([0,d]\). If
    \(f(q)[m+1]_q=\sum_i c_iq^i\), then
    \[
         c_i-c_{i-1}=b_i-b_{i-m-1}.
    \]
    For \(1\leq i\leq(d+m)/2\), we have \(i\leq d\), and
    \(\min(i,d-i)>i-m-1\). Palindromicity and strict increase up to
    the middle therefore give \(c_i>c_{i-1}\).
    For the second assertion, \(j<(d+t)/2\leq d\) and
    \(\min(j,d-j)>j-t\). The same argument, including the zero
    extension when \(j<t\), gives \(b_j>b_{j-t}\).
\end{proof}

For later use, \([n]_q!\) is strictly unimodal for \(n\geq3\).
Check \([3]_q!\) and \([4]_q!\) directly, and apply
\Cref{lem:multiply} to the remaining factors: for \(a\geq5\),
\[
(a-1)(a-2)/2>a-1.\]
Likewise
\[
 F_m(q):=\prod_{j=1}^m[2j]_q
\]
is strictly unimodal for \(m\geq3\). Indeed,
\[
 F_3(q)=1+3q+5q^2+7q^3+8q^4+8q^5
        +7q^6+5q^7+3q^8+q^9,
\]
and \(m^2>2m+1\) for \(m\geq3\), so the lemma applies
when adjoining \([2m+2]_q\).
\section{Singular loci and coefficient differences}
\label{section:singular-locus-asymmetry}

We first relate Kazhdan--Lusztig heights at maximal rationally
singular elements to ascending reflections. We then express
\(D_w(q)\) in terms of the rationally singular locus.
These give the coefficient bound in \Cref{thm:sharp-asymmetry}
and the length bound at equality in Section~\ref{subsec:sharp-bounds}.

\subsection{A lower bound at maximal rationally singular elements}

\begin{proposition}
\label{prop:bruhat-degree-identity}
    Let \(v\in\operatorname{MaxRatSing}(y)\), put \(d=d_y(v)\) and
    \(r=h_{\mathrm{KL}}(v,y)\), and write
    \[
        P_{v,y}(q)
        =
        1+\sum_{j\geq r}c_jq^j,
        \qquad
        c_j\in\mathbb Z_{\geq0},
        \qquad
        c_r>0.
    \]
    Then
    \begin{equation}
    \label{eq:bruhat-degree-coefficients}
        \deg_y^+(v)-d
        =
        \sum_{j\geq r}c_j(d-2j)>0.
    \end{equation}
    Consequently,
    \[
        h_{\mathrm{KL}}(v,y)
        \geq
        d_y(v)-\frac12\deg_y^+(v).
    \]
    Equality holds if and only if
    \[
        P_{v,y}(q)=1+q^{h_{\mathrm{KL}}(v,y)}.
    \]
\end{proposition}

\begin{proof}
    Because \(v\) is maximal in
    \(\{x\leq y:P_{x,y}(q)\neq1\}\), one has
    \[
        P_{z,y}(q)=1
        \qquad
        \text{for every }v<z\leq y.
    \]
    By \Cref{lem:KL-derivative-identity},
    \[
        dP_{v,y}(1)-2P'_{v,y}(1)
        =
        \sum_{\substack{t\in T\\v<tv\leq y}}P_{tv,y}(1)
        =
        \deg_y^+(v).
    \]
    Substituting the displayed expansion of \(P_{v,y}(q)\) gives
    \eqref{eq:bruhat-degree-coefficients}.
    
    The defining degree bound for Kazhdan--Lusztig polynomials
    \cite[\S2]{KL79} gives \(j<d/2\) whenever \(c_j>0\).
    Since \(c_r\geq1\), we obtain
    \[
         \deg_y^+(v)-d
         =\sum_{j\geq r}c_j(d-2j)\geq d-2r.
    \]
    This is the asserted lower bound. Every term with \(c_j>0\)
    is positive, so equality holds exactly when \(c_r=1\)
    and \(c_j=0\) for \(j>r\).
\end{proof}

\subsection{From singular loci to asymmetry}
\label{subsec:singular-locus-asymmetry}

For an element \(w\) of a Coxeter system, define
\[
    S_w(q)=\sum_{x\in\operatorname{RatSing}(w)}q^{\ell(x)}.
\]
This polynomial records the lengths of the elements in
\(\operatorname{RatSing}(w)\). Under the value-two hypothesis,
it determines \(D_w(q)\) through the following identity.

\begin{proposition}\label{prop:whole-locus-transfer}
    Let \((W,S)\) be an arbitrary Coxeter system and suppose that
    \(P_{\mathrm{id},w}(q)=1+q^h\), where \(h>0\).
    Then
    \begin{equation}\label{eq:whole-locus-transfer}
        D_w(q)=q^hS_w(q)-q^{\ell(w)-h}S_w(q^{-1}).
    \end{equation}
    If \(\operatorname{MaxRatSing}(w)=\{v\}\) and \(P_v(q)\)
    is palindromic, put \(d=\ell(w)-\ell(v)\). Then
    \begin{equation}\label{eq:smooth-component-transfer}
         D_w(q)=q^h(1-q^{d-2h})P_v(q).
    \end{equation}
\end{proposition}

\begin{proof}
    By \Cref{prop:monotonicity},
    \[
        1\leq P_{x,w}(q)\leq 1+q^h
        \qquad(x\leq w)
    \]
    coefficientwise. Hence \(P_{x,w}=1+q^h\) on
    \(\operatorname{RatSing}(w)\), and \(P_{x,w}=1\) elsewhere.
        Therefore
    \[
        \sum_{x\leq w}q^{\ell(x)}P_{x,w}(q)
        =P_w(q)+q^hS_w(q).
    \]
    By \cite[Lemma~2.6(v)]{KL79}, the left-hand side is
    palindromic of degree \(\ell(w)\). Hence
    \[
        P_w(q)+q^hS_w(q)
        =q^{\ell(w)}P_w(q^{-1})
         +q^{\ell(w)-h}S_w(q^{-1}),
    \]
    which proves \eqref{eq:whole-locus-transfer}.
    
    Monotonicity makes \(\operatorname{RatSing}(w)\) a lower
    order ideal. If its unique maximal element is \(v\), then
    \(\operatorname{RatSing}(w)=[\mathrm{id},v]\), so \(S_w=P_v\).
    Palindromicity gives
    \[
        q^{\ell(v)}P_v(q^{-1})=P_v(q).
    \]
    Substituting into \eqref{eq:whole-locus-transfer} gives
    \eqref{eq:smooth-component-transfer}.
\end{proof}

The coefficient of \(q^k\) in \eqref{eq:whole-locus-transfer}
is the number of rationally singular elements of length \(k-h\),
minus the number of length \(\ell(w)-k-h\).
Bounding the first count by all elements of \(W\) gives
the following coefficient bound.

\begin{corollary}\label{cor:universal-asymmetry-bound}
    Let \(W\) be a finite Weyl group with longest element
    \(w_0\), and suppose that
    \(P_{\mathrm{id},w}=1+q^h\). For \(1\leq k<\ell(w)/2\),
    \[
        0\leq[q^k]D_w(q)\leq[q^{k-h}]P_{w_0}(q)
        \leq[q^{k-1}]P_{w_0}(q),
    \]
    where coefficients in negative degrees are zero.
    For \(k\geq h\), equality in the first upper bound holds if and
    only if both of the following conditions hold:
    \begin{enumerate}[label=(\roman*)]
        \item Every element of \(W\) of length at most \(k-h\) belongs to
        \(\operatorname{RatSing}(w)\).
        \item \(d_w(v)\geq k+h+1\) for every
        \(v\in\operatorname{MaxRatSing}(w)\).
    \end{enumerate}
\end{corollary}

\begin{proof}
    Write \(S_w(q)=\sum_j c_jq^j\), and put \(c_j=0\)
outside its degree range. By
    \eqref{eq:whole-locus-transfer},
    \[
        [q^k]D_w(q)=c_{k-h}-c_{\ell(w)-k-h}
        \leq c_{k-h}\leq[q^{k-h}]P_{w_0}(q).
    \]
    Nonnegativity is \cite[Theorem~A]{BE09}. The length distribution of
    a finite Weyl group is unimodal, either by hard Lefschetz for its
    complete flag variety or by its product formula in \(q\)-integers.
    Since \(k-1<\ell(w)/2\leq\ell(w_0)/2\), this gives the last bound.
    For equality, the required two conditions are
    \(c_{k-h}=[q^{k-h}]P_{w_0}(q)\) and \(c_{\ell(w)-k-h}=0\).
    The first is equivalent to (i), because the rationally singular
    locus is a lower ideal and every element extends along a saturated
    chain to \(w_0\). The index \(\ell(w)-k-h\) is positive.
    A nonempty lower ideal contains elements of every length between
    zero and its maximum length, so its coefficient there vanishes
    exactly when every maximal element has smaller length. This is (ii).
\end{proof}

\begin{corollary}
\label{lem:noninversion-codimension}
    Let \(W\) be a finite Weyl group with longest element
    \(w_0\). If
    \(v\in\operatorname{MaxRatSing}(w)\) and
    \(h=h_{\mathrm{KL}}(v,w)\), then
    \[
         d_w(v)\leq \ell(w_0)-\ell(w)+2h.
    \]
\end{corollary}
\begin{proof}
    Put \(d=d_w(v)\). Of the \(\ell(w_0)\) reflections in \(W\), exactly
    \(\ell(v)\) decrease the length of \(v\). By
    \Cref{prop:bruhat-degree-identity},
    \[
         2d-2h\leq\deg_w^+(v)
         \leq \ell(w_0)-\ell(v)=\ell(w_0)-\ell(w)+d.
    \]
    Rearranging proves the claim.
\end{proof}
\section{Sharp asymmetry bounds and equality}
\label{section:asymmetry}

\subsection{The coefficient bound and its equality case}
\label{subsec:sharp-bounds}

If equality in \Cref{thm:sharp-asymmetry} holds at degree \(k\),
then \(P_{\mathrm{id},w}(q)=1+q\) and
\[
    \ell(w)\leq\ell(w_0)-k,
    \qquad
    N(w)\leq\left\lceil\frac{\ell(w_0)-k}{2}\right\rceil-1.
\]
The families below attain both bounds.

\begin{proof}[Proof of \Cref{thm:sharp-asymmetry} and the equality assertions above]
    Write \(P_{\mathrm{id},w}=1+q^h\).
    \Cref{cor:universal-asymmetry-bound} gives
    \[
        [q^k]D_w(q)\leq[q^{k-h}]P_{w_0}(q)\leq[q^{k-1}]P_{w_0}(q).
    \]
    Strict unimodality of \(P_{w_0}(q)\) forces \(h=1\) if the final
    bound is attained and \(k\geq2\). For \(k=1\), the transfer
    formula gives \([q]D_w(q)=0\) when \(h>1\), so equality
    again forces \(h=1\).
    
    Suppose equality holds. For any
    \(v\in\operatorname{MaxRatSing}(w)\), monotonicity gives
    \(P_{v,w}=1+q\). The equality criterion in
    \Cref{cor:universal-asymmetry-bound}, together with
    \Cref{lem:noninversion-codimension}, gives
    \[
        k+2\leq d_w(v)\leq \ell(w_0)-\ell(w)+2.
    \]
    Hence \(\ell(w)\leq \ell(w_0)-k\), and counting positive integers
    below \(\ell(w)/2\) gives the bound on \(N(w)\).
    
    For sharpness, take \(r=k\) in \Cref{thm:type-A-uniform}
    or \Cref{thm:type-B-sharp-family} below, respectively.
    The resulting element has length \(\ell(w_0)-k\), attains the first
    \(k\) coefficient bounds, and has positive differences in
    every degree below the middle. In the stated ranges,
    \(\ell(w_0)>3k\), so those first \(k\) degrees lie below the middle.
\end{proof}

\subsection{A family from one rank condition}

We construct a type-\(A\) family attaining the coefficient bound.
For \(n\geq4\) and \(1\leq r\leq n-3\), define
\[
    \begin{aligned}
        w_r^n
        &=n\,(n-2)\cdots(n-r-1)\,(n-1)
          (n-r-2)\cdots1,\\
        v_r^n
        &=(n-2)(n-3)\cdots(n-r-2)\,n\,(n-1)
          (n-r-3)\cdots1.
    \end{aligned}
\]
Empty decreasing strings are omitted. Throughout this subsection,
write \(w_r=w_r^n\) and \(v_r=v_r^n\), with \(n\) fixed.
The permutation \(w_r\) is obtained from \(w_0\) by moving
\(n-1\) past the next \(r\) entries, so
\(\ell(w_r)=\binom n2-r\).
The permutation \(v_r\) will index its singular locus.

\begin{theorem}\label{thm:type-A-uniform}
    For \(n\geq4\) and \(1\leq r\leq n-3\), the following hold.
    \begin{enumerate}[label=(\roman*)]
        \item The singular locus of \(X_{w_r}\) is the smooth
        variety \(X_{v_r}\), and
        \[
            d_{w_r}(v_r)=r+2.
        \]

        \item For every \(x\leq w_r\),
        \[
            P_{x,w_r}(q)=
            \begin{cases}
                1+q,&x\leq v_r,\\
                1,&x\not\leq v_r.
            \end{cases}
        \]

        \item The Poincar\'e and asymmetry polynomials satisfy
        \[
            P_{v_r}(q)=[n-2]_q!\,[n-r-2]_q[n-r-1]_q,
            \qquad
            D_{w_r}(q)=q(1-q^r)P_{v_r}(q).
        \]

        \item The coefficients satisfy
        \[
            [q^i]D_{w_r}(q)>0
            \qquad\left(1\leq i<\frac{\ell(w_r)}{2}\right),
        \]
        and attain the coefficient bound:
        \[
            [q^i]D_{w_r}(q)=[q^{i-1}][n]_q!
            \qquad\bigl(1\leq i\leq\min(r,n-r-2)\bigr).
        \]
    \end{enumerate}
\end{theorem}

\begin{proof}
    \noindent\textup{(i)}
    Put \(p=r+1\), fix a complete flag \(E_\bullet\) in \(\mathbb C^n\),
    and let \(\pi\) send a complete flag \(F_\bullet\) to \(F_p\).
    In \(\operatorname{Gr}(p,n)\), consider
    \[
         Y=\{U:\dim(U\cap E_{n-2})\geq p-1\},
        \qquad Y_0=\operatorname{Gr}(p,E_{n-2}).
    \]
    Schubert cells in \(\operatorname{Gr}(p,n)\) are indexed by
    increasing \(p\)-subsets of \(\{1,\ldots,n\}\), relative to
    \(E_\bullet\), with Bruhat order given by componentwise comparison. 
    The condition on \(Y\) requires at least \(p-1\)
    selected indices to be at most \(n-2\). Its largest index is
    therefore \(\{n-p,n-p+1,\ldots,n-2,n\}\).
    The largest permutation over this index is obtained by
    listing the subset decreasingly, followed by its complement
    decreasingly. This gives \(w_r\), so
    \(\pi^{-1}(Y)=X_{w_r}\).
    Applying the same procedure to the largest \(p\)-subset
    of \(\{1,\ldots,n-2\}\) gives
    \(\pi^{-1}(Y_0)=X_{v_r}\).
    
    The condition on \(Y\) says that
    \(U\to\mathbb C^n/E_{n-2}\) has rank at most one, and
    \(Y_0\) is its rank-zero locus.
    At a rank-one point, suitable bases give the quotient map
    in a Grassmannian chart the matrix
    \[
        \begin{pmatrix}
            1&0&\cdots&0\\
            z_1&z_2&\cdots&z_p
        \end{pmatrix}.
    \]
    The rank condition is \(z_2=\cdots=z_p=0\), so
    \(Y\setminus Y_0\) is smooth.
    At a point of \(Y_0\), a Grassmannian chart writes a nearby
    subspace as the graph of a linear map from \(U\) to a complement.
    Its component in \(\mathbb C^n/E_{n-2}\) is an arbitrary
    \(2\times p\) matrix; the other coordinates are free.
    Thus a transverse slice along \(Y_0\) is the cone of
    \(2\times p\) matrices of rank at most one, of dimension \(p+1\).
    At its vertex, all defining \(2\times2\) minors have zero
    differential, so the tangent space has dimension \(2p>p+1\),
    since \(p\geq2\). Thus every point of \(Y_0\) is singular, and
    \[
        \operatorname{Sing}(Y)=Y_0,
        \qquad
        \operatorname{codim}_Y Y_0=p+1=r+2.
    \]
    
    Since \(\pi\) is locally trivial with smooth fiber,
    \[
        \begin{aligned}
            \operatorname{Sing}(X_{w_r})
            &=\pi^{-1}(Y_0)=X_{v_r},\\
            d_{w_r}(v_r)
            &=\operatorname{codim}_Y Y_0=r+2.
        \end{aligned}
    \]
    Moreover, \(X_{v_r}\) is smooth, being a flag bundle over the
    Grassmannian \(Y_0\). This proves (i).
    
    \smallskip\noindent\textup{(ii)}
    Let \(\rho:\mathbb C^n\to\mathbb C^n/E_{n-2}\) be the quotient map.
    Consider the incidence variety
    \[
         \widetilde Y
         =\{(U,\lambda):U\in\operatorname{Gr}(p,n),\
               \lambda\in\mathbb P(\mathbb C^n/E_{n-2}),\
               \rho(U)\subseteq\lambda\}.
    \]
    For fixed \(\lambda\), the space \(\rho^{-1}(\lambda)\) has
    dimension \(n-1\). Thus \(\widetilde Y\) is a smooth
    \(\operatorname{Gr}(p,n-1)\)-bundle over \(\mathbb P^1\).
    The map \(\widetilde Y\to Y\) is proper and is an isomorphism
    off \(Y_0\), where \(\lambda=\rho(U)\).
    Its fiber over \(Y_0\) is \(\mathbb P^1\).
    
    A resolution is small if the locus of fibers of dimension at
    least \(a\) has codimension greater than \(2a\) for every \(a>0\).
    Here the only positive fiber dimension is one, and
    \(\operatorname{codim}_Y Y_0=r+2>2\).
    Pulling back along \(\pi\) gives a small resolution of \(X_{w_r}\)
    whose fibers are \(\mathbb P^1\) over points of \(X_{v_r}\)
    and points elsewhere. The small-resolution theorem and the
    Kazhdan--Lusztig interpretation of local intersection cohomology
    \cite{GM83,KL80} therefore give \(P_{x,w_r}(q)=1+q\)
    for \(x\leq v_r\), and \(P_{x,w_r}(q)=1\) otherwise.
    This proves (ii).
    
    \smallskip\noindent\textup{(iii)}
    Over \(U\in Y_0\), the fiber of \(\pi\) is the product of the
    complete flag varieties of \(U\) and \(\mathbb C^n/U\).
    Their Schubert cell decompositions give
    \[
         P_{v_r}(q)
         =\genfrac{[}{]}{0pt}{}{n-2}{r+1}_q
               [r+1]_q![n-r-1]_q!
         =[n-2]_q![n-r-2]_q[n-r-1]_q.
    \]
    
    Now \eqref{eq:smooth-component-transfer}, with \(h=1\) and
    \(d=r+2\), gives the claimed formula for \(D_{w_r}\).
    
    \smallskip\noindent\textup{(iv)}
    Put \(f=P_{v_r}\) and \(e=\deg f=\binom n2-2r-2\).
    The polynomial \(f\) is strictly unimodal. For \(n=4\), it is
    \((1+q)^2\). For \(n=5\), its coefficient lists are
    \((1,4,8,10,8,4,1)\) and \((1,3,4,3,1)\).
    For \(n\geq6\), apply \Cref{lem:multiply} to \([n-2]_q!\)
    and the two remaining factors, omitting any factor \([1]_q\).
    Their degrees are at most \(n-3<(n-2)(n-3)/2\).
    Moreover,
    \[
         e-r\geq\binom n2-3(n-3)-2
               =\frac{n^2-7n+14}{2}>0\qquad(n\geq4).
    \]
    Since \(D_{w_r}=q(1-q^r)f\), the second assertion of
    \Cref{lem:multiply} gives positivity for
    \(1\leq i<(e+r+2)/2=\ell(w_r)/2\).
    
    Finally, put \(k=\min(r,n-r-2)\).
    The factor \(1-q^r\) does not affect degrees below \(k\), so
    \(D_{w_r}\equiv qP_{v_r}\pmod{q^{k+1}}\).
    Each of \(n-r-2,n-r-1,n-1,n\) is at least \(k\).
    Their \(q\)-integers therefore agree modulo \(q^k\), giving
    \[
         P_{v_r}(q)\equiv[n]_q!\pmod{q^k}.
    \]
    Comparing coefficients proves the last assertion.
\end{proof}

\subsection{The same rank condition in type \texorpdfstring{\(B\)}{B}}
\label{subsec:bottom-asymmetry}

The rank condition used in type \(A\) also gives sharp bounds
in type \(B\). We realize the common Coxeter system of types
\(B_n\) and \(C_n\) using symplectic flags. Length, Bruhat order, and Kazhdan--Lusztig polynomials depend
only on this Coxeter system, so the resulting formulas for
\(P_w\), \(D_w\), and \(P_{x,w}\) apply in both types. The resolution again
has fibers that are points or projective lines. The rationally
singular locus need not be rationally smooth, so we use the
whole-locus transfer formula.

For \(n\geq3\) and \(1\leq r\leq n-2\), define
\[
    c_r^n=\overline1\cdots\overline{n-r-2}\,
          \overline{n-1}\,
          \overline{n-r-1}\cdots\overline{n-2}\,\overline n.
\]
Empty strings are omitted. Write \(c_r=c_r^n\) throughout this
subsection, with \(n\) fixed. This permutation has length
\(n^2-r\); the proof below verifies this and computes its
asymmetry polynomial.

The two families arise from the same element of the standard
parabolic subgroup \(S_n\subset W(B_n)\).
For this comparison, write \(w_0^A\) and \(w_0^B\) for the
longest elements of \(S_n\) and \(W(B_n)\), respectively.
For \(1\leq r\leq n-3\), put
\[
    t_r=s_{n-2}s_{n-3}\cdots s_{n-r-1}.
\]
Then
\[
    w_r=t_rw_0^A,
    \qquad c_r=t_rw_0^B.
\]
Since \(\ell(t_r)=r\), both varieties have codimension \(r\)
in their respective flag varieties.
Multiplication by the longest element identifies their upper
Bruhat intervals with the order dual of \([\mathrm{id},t_r]\).
Their lower intervals require separate calculations, which
we carry out geometrically.

\begin{theorem}\label{thm:type-B-sharp-family}
    For \(n\geq3\) and \(1\leq r\leq n-2\), we have
    \(P_{\mathrm{id},c_r}(q)=1+q\) and
    \[
         D_{c_r}(q)
         =q(1-q^r)F_{n-2}(q)[2n-r-3]_q[2n-r-1]_q.
    \]
    Moreover,
    \[
         [q^i]D_{c_r}(q)>0
         \qquad\left(1\leq i<\frac{n^2-r}{2}\right),
    \]
    and
    \[
         [q^i]D_{c_r}(q)=[q^{i-1}]F_n(q)
         \qquad(1\leq i\leq r).
    \]
\end{theorem}

We fix the symplectic space and coordinate conventions used in this proof and in Section~\ref{sec:type-B-families}.
Let \(V\) be a complex vector space with basis
\[e_1,\ldots,e_n,f_1,\ldots,f_n\] and alternating form
\(\langle e_i,f_j\rangle=\delta_{ij}\), with zero pairings
among the \(e_i\)'s and among the \(f_i\)'s.
For a subspace \(U\subseteq V\), write
\[
    U^\perp=\{v\in V:\langle u,v\rangle=0
    \text{ for every }u\in U\}.
\]
The subspace \(U\) is isotropic if \(U\subseteq U^\perp\),
and the radical of the restricted form is \(U\cap U^\perp\).
Write \(\operatorname{IG}(k,V)\) for the variety of isotropic
\(k\)-planes in \(V\).
Set \(E_0=0\), \(E_i=\langle e_1,\ldots,e_i\rangle\)
for \(1\leq i\leq n\), and
\(E_i=E_{2n-i}^\perp\) for \(n<i\leq2n\).
Thus \(E_\bullet\) is the coordinate flag for the ordered basis
\(e_1,\ldots,e_n,f_n,\ldots,f_1\).
Put \(G=\operatorname{Sp}(V)\), and let \(B\) stabilize
\(E_\bullet\). Thus \(G/B\) parametrizes complete isotropic flags in \(V\).
The fixed flag of a signed permutation \(u\) has its
\(i\)-plane spanned by
\(\varepsilon_1(u),\ldots,\varepsilon_i(u)\), for
\(1\leq i\leq n\), where
\begin{equation}\label{eq:signed-fixed-flag}
    \varepsilon_i(u)=
    \begin{cases}
        e_{n+1-u(n+1-i)},&u(n+1-i)>0,\\
        f_{n+1-|u(n+1-i)|},&u(n+1-i)<0.
    \end{cases}
\end{equation}

\begin{proof}
    Put \(p=r+1\), \(m=2n-p\), and \(w=c_r^{-1}\).
    Inversion preserves length and Bruhat order, so
    \(D_w=D_{c_r}\). It also preserves Kazhdan--Lusztig
    polynomials, so it suffices to prove the assertions for \(w\).
    Since \(w=w_0 s_{n-p}\cdots s_{n-2}\), we have
    \[
        \ell(w)=n^2-p+1=n^2-r.
    \]

    \smallskip\noindent\emph{Computing the Kazhdan--Lusztig polynomials.}
    Set \(H=E_m\), whose radical is \(R=E_p\).
    In the Grassmannian of isotropic two-planes, consider
    \[
        Y=\{U\in\operatorname{IG}(2,V):U\cap H\neq0\},
        \qquad
        Y_0=\{U\in\operatorname{IG}(2,V):U\subset H\}.
    \]
    The fixed two-plane of \(w\) is
    \(\langle f_1,f_{p+1}\rangle\), which belongs to
    \(Y\setminus Y_0\).
    Choosing a line in \(U\cap H\) gives
    \[
        Z=\{(\lambda,U)\in\mathbb P(H)\times
                    \operatorname{IG}(2,V):\lambda\subset U\}
        \longrightarrow Y.
    \]
    Over \(\lambda\in\mathbb P(H)\), the choices of \(U\) form
    \(\mathbb P(\lambda^\perp/\lambda)\simeq\mathbb P^{2n-3}\).
    Thus \(Z\) is a smooth irreducible projective bundle.
    The map \(Z\to Y\) is proper, is an isomorphism off \(Y_0\),
    and has fiber \(\mathbb P^1\) over \(Y_0\).
    Consequently \(Y\) is irreducible of dimension \(m+2n-4\).

    The projection \(\pi\colon G/B\to\operatorname{IG}(2,V)\)
    is locally trivial. Over \(U\), its fiber is the product
    of \(\mathbb P(U)\) and the complete symplectic flag
    variety of \(U^\perp/U\). This fiber is smooth and
    irreducible, of dimension \(1+(n-2)^2\) and
    Poincar\'e polynomial \([2]_qF_{n-2}(q)\).
    Hence \(\pi^{-1}(Y)\) is irreducible of dimension
    \[
        m+2n-4+1+(n-2)^2=\ell(w).
    \]
    Since \(Y\) is closed and \(B\)-stable and contains the
    fixed two-plane of \(w\), we have
    \(X_w\subseteq\pi^{-1}(Y)\); equality follows by dimension.

    The alternating form on \(H\) has rank \(2(n-p)>0\).
    Under the Pl\"ucker embedding
    \(\operatorname{span}\{u_1,u_2\}\mapsto[u_1\wedge u_2]\),
    isotropy is the nonzero linear equation
    \(\langle u_1,u_2\rangle=0\).
    Thus \(Y_0\) is a hyperplane section of
    \(\operatorname{Gr}(2,H)\), so \(\dim Y_0=2m-5\) and
    \[
        \operatorname{codim}_Y Y_0=p+1>2.
    \]
    Therefore \(Z\to Y\), and its smooth pullback to \(X_w\),
    are small resolutions. Their fibers give
    \(P_{x,w}=1+q\) on \(\pi^{-1}(Y_0)\) and \(1\) elsewhere
    \cite{GM83,KL80}.
    Since \(E_2\subset H\), in particular
    \(P_{\mathrm{id},w}=1+q\).

    \smallskip\noindent\emph{Computing the asymmetry polynomial.}
    Count the flags \(\lambda\subset U\subset H\), with
    \(\lambda\) a line and \(U\) an isotropic two-plane.
    Over \(U\in Y_0\), the fiber is \(\mathbb P^1\).
    Over \(\lambda\in\mathbb P(H)\), the fiber is
    \(\mathbb P((H\cap\lambda^\perp)/\lambda)\):
    it is \(\mathbb P^{m-2}\) if \(\lambda\subset R\), and
    \(\mathbb P^{m-3}\) otherwise.
    Write \(P_{Y_0}(q)=\sum_{C\subseteq Y_0}q^{\dim C}\),
    where \(C\) runs over the Schubert cells contained in \(Y_0\).
    Since \(Y_0\) is closed and \(B\)-stable, this is its
    ordinary Poincar\'e polynomial. Counting over finite fields gives
        \[
        \begin{aligned}
            [2]_qP_{Y_0}(q)
            &=[p]_q[m-1]_q+([m]_q-[p]_q)[m-2]_q\\
            &=[m]_q[m-2]_q+q^{m-2}[p]_q.
        \end{aligned}
    \]
    
    Over each Schubert cell of \(Y_0\), the projection \(\pi\)
    is a trivial bundle with the fiber described above.
    Since \(\pi^{-1}(Y_0)\) is the rationally singular locus
    of \(X_w\), counting its cells gives
    \[
        \begin{aligned}
            S_w(q)
            &=[2]_qF_{n-2}(q)P_{Y_0}(q)\\
            &=f(q)+q^{m-2}F_{n-2}(q)[p]_q,
        \end{aligned}
        \qquad
        f(q):=F_{n-2}(q)[m-2]_q[m]_q.
    \]
    The polynomial \(f\) is palindromic of degree \(n^2-2p\).
    Since
    \[
        2(m-1)+\deg(F_{n-2}(q)[p]_q)=\ell(w),
    \]
    the second term of \(qS_w(q)\) is symmetric about \(\ell(w)/2\)
    and cancels in \eqref{eq:whole-locus-transfer}. Hence
    \[
        \begin{aligned}
            D_w(q)
            &=qS_w(q)-q^{\ell(w)-1}S_w(q^{-1})\\
            &=(q-q^p)f(q)=q(1-q^{p-1})f(q).
        \end{aligned}
    \]
    Substituting \(p=r+1\) proves the asserted factorization.

    \smallskip\noindent\emph{The coefficient claims.}
    Put \(a=m-2\geq n-1\).
    The polynomial \([2]_q[a]_q[a+2]_q\) has degree \(2a+1\)
    and coefficients
    \(1,3,5,\ldots,2a-1,2a\) up to the middle.
    It is therefore strictly unimodal.
    Adjoining the remaining factors \([2j]_q\),
    \(2\leq j\leq n-2\), preserves strict unimodality by
    \Cref{lem:multiply}, since
    \(2j-1\leq2n-5<2a+1\).
    Hence \(f\) is strictly unimodal.
    As \(r\leq\deg f=n^2-2r-2\), the same lemma applied to
    \(D_w=q(1-q^r)f\) gives
    \[
        [q^i]D_w(q)>0
        \qquad\left(1\leq i<\frac{n^2-r}{2}\right).
    \]
    Finally, \(m-2,m,2n-2,2n\geq r\), so
    \[
        f\equiv F_n\pmod{q^r},
        \qquad
        D_w\equiv qF_n\pmod{q^{r+1}}.
    \]
    Comparing coefficients proves the remaining equalities.
\end{proof}

\begin{proof}[Proof of \Cref{thm:number-asymmetries}]
    The longest element has no coefficient differences.
    Every other element has length at most \(\ell(w_0)-1\), so
    \[
        N(w)\leq\left\lceil\frac{\ell(w_0)-1}{2}\right\rceil-1
        =\left\lfloor\frac{\ell(w_0)}2\right\rfloor-1.
    \]
    By definition, \(u_n=w_1^n\) in type \(A\) and \(u_n=c_1^n\)
    in type \(B\). By \Cref{thm:type-A-uniform,thm:type-B-sharp-family},
    \(P_{\mathrm{id},u_n}(q)=1+q\), \(\ell(u_n)=\ell(w_0)-1\),
    and every coefficient of \(D_{u_n}(q)\) in degrees
    \(1\leq i<\ell(u_n)/2\) is positive. Thus \(u_n\) attains the bound.
\end{proof}
\section{Realizing \texorpdfstring{\(P_{\mathrm{id},w}(q)=1+q^h\)}{Pid,w(q)=1+qh} in type \texorpdfstring{\(B\)}{B}}
\label{sec:type-B-families}

This section proves \Cref{thm:type-B-height-spectrum}.
The families in \Cref{thm:type-A-uniform,thm:type-B-sharp-family}
satisfy \(P_{\mathrm{id},w}(q)=1+q\).
We now construct elements
in \(W(B_n)\) such that \(P_{\mathrm{id},w}(q)=1+q^h\) for every
\(1\leq h\leq2n-4\).

We use symplectic flags \(L\subset M\) of dimensions one and two.
The incidence conditions below produce a transverse quadratic
cone of dimension \(2h+1\) along the singular locus.
Its normalized local intersection-cohomology polynomial is
\(1+q^h\), which gives the desired Kazhdan--Lusztig polynomial.
We first specify the signed permutations indexing the resulting
Schubert variety and its rationally singular locus.
Let \(n\geq3\) and \(1\leq h\leq2n-4\), and define
\[
    t=
    \begin{cases}
        n-h-1,&h\leq n-2,\\
        n-h-2,&h\geq n-1,
    \end{cases}
    \qquad
    s=
    \begin{cases}
        h-n+2,&h\leq n-3,\\
        h-n+3,&h\geq n-2.
    \end{cases}
\]
Let \(i_1<\cdots<i_{n-2}\) and \(j_1<\cdots<j_{n-2}\)
list the elements of
\[
    \{1,\ldots,n\}\setminus\{|t|,n-1\},
    \qquad
    \{1,\ldots,n\}\setminus\{|s|,n\},
\]
respectively. Define the signed permutations
\begin{equation}\label{eq:type-B-height-permutations}
    b_h^n=\overline{i_1}\cdots\overline{i_{n-2}}\,
          \overline{n-1}\,t,
    \qquad
    a_h^n=\overline{j_1}\cdots\overline{j_{n-2}}\,s\,n.
\end{equation}
For fixed \(n\), write \(b_h=b_h^n\) and \(a_h=a_h^n\).
For example, \(n=5\) and \(h=2\) give
\[
 b_2^5=\overline1 \overline3 \overline5 \overline4 2,
 \qquad a_2^5=\overline2 \overline3 \overline4 \overline1 5.
\]

The proof also gives
\[
    \ell(b_h)=(n-1)^2+h,
    \qquad \ell(b_h)-\ell(a_h)=2h+1.
\]
For \(n\geq5\), the coefficients of \(D_{b_h}\) vanish before
degree \(h\) and are positive in every degree from \(h\) to
strictly below \(\ell(b_h)/2\).
\begin{proof}[Proof of \Cref{thm:type-B-height-spectrum}]
    We use the symplectic setup from
    Section~\ref{subsec:bottom-asymmetry}.
    First we realize \(X_{b_h}\) as a flag bundle over an incidence
    variety \(Z\). We then compute its Kazhdan--Lusztig polynomials
    and identify \(X_{a_h}\) as its rationally singular locus.
    Finally, we compute \(D_{b_h}\) and establish the positivity
    assertion above.
    
    \smallskip\noindent\emph{1. Realizing \(X_{b_h}\) as a flag bundle.}
    Put
    \[
        K=E_{h+2},\qquad H=E_1^\perp.
    \]
    The range of \(h\) ensures \(E_1\subset K\subset H\).
    In the variety of isotropic flags of dimensions one and two, set
    \[
        Z=\{(L,M):L\subset M,\ L\subset K,\ M\subset H,\
                    M\text{ isotropic}\}.
    \]
    The projection to \(L\) has fiber
    \(\mathbb P((H\cap L^\perp)/L)\).
    Since \(H^\perp=E_1\), pairing with a nonzero vector in
    \(L\) gives a nonzero functional on \(H\) when \(L\ne E_1\).
    Thus the fiber is \(\mathbb P^{2n-4}\) for \(L\ne E_1\)
    and \(\mathbb P^{2n-3}\) for \(L=E_1\).
    Thus the part over \(\mathbb P(K)\setminus\{E_1\}\) is an
    irreducible projective bundle of dimension \(2n+h-3\).
    
    Before imposing isotropy, the same flags form a smooth projective
    bundle with fiber \(\mathbb P(H/L)\).
    As a hypersurface in this smooth bundle, \(Z\) has pure
    dimension \(2n+h-3\). The fiber over \(E_1\) has dimension
    \(2n-3<2n+h-3\), so every irreducible component meets
    the irreducible open subset over \(L\ne E_1\).
    Hence \(Z\) is irreducible.
    It is closed and \(B\)-stable; since the partial flag variety
    has finitely many \(B\)-orbits, \(Z\) is a Schubert variety.
    
    We now pass to complete flags and identify the inverse image of \(Z\).
    Let \(J=\{s_0,\ldots,s_{n-3}\}\), and let \(P_J\) be the
    corresponding standard parabolic subgroup of \(G\).
    The projection
    \[
        \pi\colon G/B\longrightarrow G/P_J
    \]
    retains the first two subspaces \(L\subset M\) of a complete
    isotropic flag, with \(\dim L=1\) and \(\dim M=2\).
    For fixed \(L\subset M\), completing the flag amounts to choosing
    a complete isotropic flag in the symplectic space \(M^\perp/M\),
    which has dimension \(2n-4\).
    Thus the fiber has dimension \((n-2)^2\) and Poincar\'e
    polynomial \(F_{n-2}\). Consequently,
    \[
        \dim\pi^{-1}(Z)
        =2n+h-3+(n-2)^2
        =(n-1)^2+h.
    \]
    
    By \eqref{eq:signed-fixed-flag}, the last entry \(t\) of
    \(b_h\) gives the last coordinate line of \(K\), and
    the preceding entry \(-(n-1)\) gives \(f_2\).
    The resulting flag satisfies the conditions of \(Z\).
    The signed length formula gives
    \[
         \ell(b_h)=
         \begin{cases}
         n^2-n-t,&t>0,\\
         n^2-n+|t|-1,&t<0,
         \end{cases}
         =(n-1)^2+h.
    \]
    Thus \(X_{b_h}\subseteq\pi^{-1}(Z)\) has the same dimension
    as this irreducible variety, proving equality.
    
    \smallskip\noindent \emph{2. Computing the Kazhdan--Lusztig polynomials.}
    The projective bundle over \(L\ne E_1\) is smooth.
    At \((E_1,\langle e_1,m\rangle)\), varying \(L\) in a
    direction \(u\in K/E_1\) changes the isotropy equation by
    \(\langle u,m\rangle\).
    The derivatives obtained by varying \(M\) with \(L=E_1\)
    vanish because \(E_1=H^\perp\).
    Thus the point is singular exactly when \(m\in K^\perp\).
    Hence
    \[
         \operatorname{Sing}(Z)
         =\{(E_1,M):E_1\subset M\subset K^\perp\}
         \simeq\mathbb P(K^\perp/E_1)
         \simeq\mathbb P^{2n-h-4}.
    \]
    Here \(E_1\subset K^\perp\subset H\), as follows from
    \(E_1\subset K\subset H\).
    
    Fix a point \((E_1,\langle e_1,m\rangle)\) of
    \(\operatorname{Sing}(Z)\).
    To determine its transverse singularity, choose a complement
    \(U\) of \(E_1\) in \(K\) and a complement \(V_0\) of
    \(\langle e_1,m\rangle\) in \(H\).
    A neighborhood of the point has coordinates
    \[
         L=\langle e_1+u\rangle,
         \qquad M=\langle e_1+u,m+z\rangle,
         \qquad u\in U,\ z\in V_0.
    \]
    Because \(m\in K^\perp\) and \(z\in H\), isotropy becomes
    \(\langle u,z\rangle=0\).
    The pairing \(U\times V_0\) has rank \(h+1=\dim U\):
    an element of \(U\) annihilating \(V_0\) also annihilates
    \(e_1\) and \(m\), hence all of \(H\); it lies in
    \(H^\perp=E_1\), so it is zero.
    After linear changes of coordinates, the equation is therefore
    \[
     u_1z_1+\cdots+u_{h+1}z_{h+1}=0,
    \]
    with \(2n-h-4\) free coordinates.
    
    The transverse cone has dimension \(2h+1\) and is the affine
    cone over a smooth projective quadric \(Q\) of dimension \(2h\).
    Write \(\beta_j=\dim H^{2j}(Q,\mathbb Q)\) and \(\beta_{-1}=0\).
    The cone formula gives the Poincar\'e polynomial of the
    normalized local intersection cohomology at the vertex:
    \[
     \sum_{j=0}^{h}(\beta_j-\beta_{j-1})q^j=1+q^h.
    \]
    The subtraction removes the image of multiplication by the
    hyperplane class. Here \(\beta_j=1\) for \(0\leq j<h\) and
    \(\beta_h=2\); see \cite{GM83} and
    \cite[Corollary, case~2]{Man01}.
    Since \(\pi\) is locally trivial with smooth fiber,
    pullback along \(\pi\) preserves these normalized local
    intersection-cohomology groups.
    Consequently, \(P_{x,b_h}=1+q^h\) when the cell \(X_x^\circ\)
    lies over \(\operatorname{Sing}(Z)\), and \(P_{x,b_h}=1\)
    otherwise.
    
    This inverse image is smooth and irreducible. The last entry of
    \(a_h\) gives \(L=E_1\), and the preceding entry gives the
    last coordinate line of \(K^\perp\). Its fixed flag therefore
    lies in the inverse image. The signed length formula gives
    \[
        \ell(a_h)=(n-2)^2+2n-h-4,
    \]
    which is the dimension of that inverse image. Thus it is
    \(X_{a_h}\), and its codimension is \(2h+1\).
    The calculation of \(P_{x,b_h}\) now gives
    \[
     \operatorname{RatSing}(b_h)=[\mathrm{id},a_h].
    \]
    Smoothness of \(X_{a_h}\) in this type-\(C\) realization
    implies rational smoothness in both types.
    
    \smallskip\noindent\emph{3. Computing the asymmetry polynomial and its coefficients.}
    We use the two fibers of \(Z\to\mathbb P(K)\) described above.
    Over a finite field with \(q\) elements, a projective space of dimension \(j\) has \([j+1]_q\) points. 
    Schubert cell counts equal the corresponding Poincar\'e polynomials. 
    Including the flag fiber gives
    \[
         \begin{aligned}
             P_{b_h}(q)
             &=F_{n-2}(q)
               \bigl(([h+2]_q-1)[2n-3]_q+[2n-2]_q\bigr)\\
             &=F_{n-2}(q)
               \bigl([h+2]_q[2n-3]_q+q^{2n-3}\bigr).
         \end{aligned}
    \]
    Since \(\operatorname{Sing}(Z)\simeq\mathbb P^{2n-h-4}\)
    and \(X_{a_h}=\pi^{-1}(\operatorname{Sing}(Z))\), the same
    flag fiber gives
    \[
        P_{a_h}(q)=[2n-h-3]_qF_{n-2}(q).
    \]
    Using \eqref{eq:smooth-component-transfer} with \(d=2h+1\),
    we obtain
    \[
         D_{b_h}=q^h(1-q)P_{a_h}
                   =(q^h-q^{2n-3})F_{n-2}.
    \]
    
    It remains to prove the positivity assertion for \(n\geq5\).
    Put \(m=n-2\) and \(c=2n-h-3\). Then
    \[
         D_{b_h}=q^h(1-q^c)F_m,
         \qquad \ell(b_h)=m^2+2h+c.
    \]
    For \(n\geq5\), the polynomial \(F_m\) is strictly unimodal,
    and \(1\leq c\leq2n-4\leq m^2=\deg F_m\).
    If \(h\leq i<\ell(b_h)/2\), then
    \(0\leq i-h<(m^2+c)/2\).
    The second assertion of \Cref{lem:multiply} therefore gives
    \([q^i]D_{b_h}(q)>0\).
    The factor \(q^h\) gives \([q^i]D_{b_h}(q)=0\) for \(i<h\).
\end{proof}

\begin{remark}
    The rank restriction in the strictness assertion matters.
    For \(n=4,h=4\), the displayed formula gives
    \([q^6]D_{b_4^4}(q)=0\), although
    \(4<6<\ell(b_4^4)/2\).
\end{remark}

Taking \(h=2n-4\) gives the following specialization.

\begin{corollary}\label{thm:type-B-height-family}
    For \(n\geq3\), write \(b=b_{2n-4}^n\) and \(a=a_{2n-4}^n\).
    Explicitly,
    \[
         b=\overline1\cdots\overline{n-3}\,
               \sbar n\,\overline{n-1}\,\overline{n-2},
         \qquad
         a=\overline1\cdots\overline{n-2}\,(n-1)n.
    \]
    Then \(\operatorname{MaxRatSing}(b)=\{a\}\),
    \(d_b(a)=4n-7\), \(\deg^+_b(a)=4n-6\), and,
    for \(x\leq b\), \(P_{x,b}=1+q^{2n-4}\) exactly when \(x\leq a\).
\end{corollary}

\begin{proof}
    Take \(h=2n-4\) in \Cref{thm:type-B-height-spectrum};
    the upward degree follows from
    \Cref{prop:bruhat-degree-identity}.
\end{proof}


\begin{thebibliography}{99}




\bibitem{BB05}
A.~Bj{\"o}rner and F.~Brenti,
\emph{Combinatorics of Coxeter groups},
Graduate Texts in Mathematics, vol.~231,
Springer, New York, 2005.

\bibitem{BL00}
S.~C.~Billey and V.~Lakshmibai,
\emph{Singular Loci of Schubert Varieties},
Progress in Mathematics, vol.~182,
Birkh{\"a}user Boston, Boston, MA, 2000.


\bibitem{BE09}
A.~Bj{\"o}rner and T.~Ekedahl,
\emph{On the shape of {B}ruhat intervals},
Ann. of Math. (2) \textbf{170} (2009), no.~2, 799--817.




\bibitem{GG24}
C.~Gaetz and Y.~Gao,
\emph{On the minimal power of \(q\) in a Kazhdan--Lusztig polynomial},
Adv. Math. \textbf{457} (2024), Paper No.~109941, 25~pp.

\bibitem{GM83}
M.~Goresky and R.~MacPherson,
\emph{Intersection homology. II},
Invent. Math. \textbf{72} (1983), no.~1, 77--129.

\bibitem{Ir88}
R.~S. Irving,
\emph{The socle filtration of a Verma module},
Ann. Sci. \'Ecole Norm. Sup. (4) \textbf{21} (1988), no.~1, 47--65.

\bibitem{KL79}
D.~Kazhdan and G.~Lusztig,
\emph{Representations of Coxeter groups and Hecke algebras},
Invent. Math. \textbf{53} (1979), no.~2, 165--184.

\bibitem{KL80}
D.~Kazhdan and G.~Lusztig,
\emph{Schubert varieties and Poincar\'e duality},
in Geometry of the Laplace Operator,
Proc. Sympos. Pure Math. \textbf{36} (1980), 185--203.

\bibitem{Kob13}
M.~Kobayashi,
\emph{Inequalities on Bruhat graphs, \(R\)- and Kazhdan--Lusztig polynomials},
J. Combin. Theory Ser. A \textbf{120} (2013), no.~2, 470--482.



\bibitem{Man01}
L.~Manivel,
\emph{Generic singularities of Schubert varieties},
arXiv:math/0105239 (2001).


\bibitem{Pla17}
D.~Plaza,
\emph{Graded cellularity and the monotonicity conjecture},
J. Algebra \textbf{473} (2017), 324--351.

\bibitem{Woo09}
A.~Woo,
\emph{Permutations with Kazhdan--Lusztig polynomial \(P_{\mathrm{id},w}(q)=1+q^h\)},
Electron. J. Combin. \textbf{16} (2009), no.~2, R10, 32~pp.
With an appendix by S.~C.~Billey and J.~Weed.



\end{thebibliography}
\end{document}